\documentclass{article}
\usepackage{amscd}
\usepackage{amssymb}

\newcommand{\md}{{\rm{Mod}}}
\newcommand{\ab}{{\rm{Ab}}}
\newcommand{\mor}{{\rm{Mor}}}
\newcommand{\ob}{{\rm{Ob}}}
\newcommand{\coker}{{\rm{coker}}}
\newcommand{\im}{{\rm{Im}}}
\newcommand{\hm}{{\rm{Hom}}}
\newcommand{\ext}{{\rm{Ext}}}
\newcommand{\set}{{\rm{Set}}}
\newcommand{\dom}{{\rm{dom}}}
\newcommand{\cod}{{\rm{cod}}}

\newtheorem{theorem}{\bf Theorem}[section]
\newtheorem{lemma}[theorem]{\bf Lemma}
\newtheorem{proposition}[theorem]{\bf Proposition}
\newtheorem{corollary}[theorem]{\bf Corollary}

\title{On Ext in the Category of Functors to \\ Preabelian Category
\footnote{In the preparation of this paper, we have been assisted 
by a grant from the 
T\"UB\.ITAK and NATO. The first author is supported by Russian Federation 
Education Department}
}
\date{}

\begin{document}
\maketitle

\begin{center}
{A. A. Husainov}\\
Komsomolsk-on-Amur State Technical University,\\
Faculty of Computer Technologies,\\
681013, Komsomolsk-on-Amur, prosp. Lenina, 27, Russia\\
husainov@knastu.ru
\end{center}
\begin{center}
{A. Pancar}\\
Ondokuz Mayis University, Faculty of Art and Sciences,\\
Department of Mathematics, 55139, Kurupelit, Samsun, Turkey\\
apancar@omu.edu.tr
\end{center}
\begin{center}
{M. Yapici}\\
Ondokuz Mayis University, Faculty of Art and Sciences,\\
Department of Mathematics, 55139, Kurupelit, Samsun, Turkey\\
myapici@omu.edu.tr
\end{center}

\smallskip
\noindent
1991 AMS Subject Classification: 18G25

\begin{abstract}
The article is devoted to extension groups in the category of functors
from a small category to an additive category with an Abelian structure
in the sense of Heller.
It is proved that
under additional assumptions
there exists a spectral sequence which converges to the extension groups.
\end{abstract}
\noindent
{Key words : \em Extension groups, cohomology of categories.}

\section{INTRODUCTION}
Let $\ab$ be the category of abelian groups, $\md_R$ the category of left
$R$-modules over a ring with identity.

Jibladze and Pirashvili \cite{jib1986}, \cite{jib1991} proved that if ${\cal A}=\md_R$
then for every
small category $\Bbb C$ and functors $F,G:{\Bbb C}\to {\cal A}$ there exists
a first quadrant spectral sequence with
$$
E^{p,q}_2=H^p({\Bbb C},\{\ext^q(F(a),G(b))\})
$$
which converges to $\ext^n(F,G)$. Here $H^n({\Bbb C},B)$ are
the $n$-th Hochschild-Mitchell cohomology groups
with coefficients in a bimodule
$B:{\Bbb C}^{op}\times {\Bbb C}\to \ab$ in the sense of \cite{bau1985},
and $\{\ext^q(F(a),G(b))\}$ is the functor
$\ext^q(F(-),G(=)): {\Bbb C}^{op}\times {\Bbb C}\to \ab$.
This spectral sequence was
generalized by Khusainov \cite{X19921} to an Abelian category $\cal A$ with
a proper class  \cite{mac1975}. We generalize it to additive categories
with Abelian structures in the sense of Heller \cite{hel1958}.

\section{ABELIAN STRUCTURE}
Let $\cal A$ be a preadditive category, $B\in \cal A$ an object.
A {\em direct sum decomposition} $(A,A',i,p,i',p')$ of $B$ consist
of objects and morphisms
$$
A \stackrel{i}\rightarrow B \stackrel{i'}\leftarrow{A'}, \quad
A \stackrel{p}\leftarrow B \stackrel{p'}\rightarrow{A'}
$$
satisfying to the conditions:
$$\begin{array}{c}
p\circ i=1_A, \quad p'\circ i'=1_{A'},\\
p'\circ i=0_{AA'},\quad p\circ i'=0_{A'A},\\
i\circ p + i'\circ p'=1_B
\end{array}$$

A  morphism $\rho:B\to A$ is called a {\em retraction} if there exists a
morphism $\nu:A\to B$ such that $\rho\circ \nu=1_A$.

An {\em additive} category is a preadditive category with a zero-object and
with finite products. We follow to Mac Lane \cite{mac1975} for the notions 
of the kernels and cokernels. 
\begin{lemma}\label{l21}
Let $\cal A$ be an additive category. Then the following conditions are
equivalent:
\begin{itemize}
\item[(1)] Every retraction has a kernel;
\item[(2)] For every morphisms $\rho:B\to A$ and $i:A\to B$ in $\cal A$
satisfying $\rho \circ i=1_A$ there are an object $A'$ and morphisms
$\rho':B\to A'$,
$i':A'\to B$ such that $(A,A',i,\rho,i',\rho')$ is a direct sum decomposition
of B.
\end{itemize}
\end{lemma}

The condition (2) is called  Cancellation Axiom in \cite{hel1958}.
Thus, if
$\cal A$ is an additive category in which all retractions have kernels
then $\cal A$ is a category with Cancellation Axiom.

The sequence of morphisms
$$\begin{CD}
0 @>>> A' @>\alpha'>> A @>\alpha''>> A'' @>>> 0
\end{CD}$$
is called a {\em short exact sequence}, abbreviated "s.e.s", if
$(A'', \alpha'')$ is a cokernel of $\alpha'$
and $(A', \alpha')$ is a kernel of $\alpha''$. More generally,
a sequence
$$
\begin{CD}
A_n@>f_n>> A_{n-1} @>f_{n-1}>>\dots @>f_1>> A_0
\end{CD}
$$
is {\em exact} if there are s.e.s
$0 \to B_j \stackrel{u_j}\to A_j \stackrel{\nu_j}\to B_{j-1} \to 0$ for
$j=1,2,\dots,n-1$, an epimorphism $\nu_n:A_n\to B_{n-1}$ and a monomorphism
$u_0:B_0\to A_0$ such that $f_j=u_{j-1}\circ \nu_j,\quad j=1,2,\dots,n$.

Let $\cal A$ be an additive category with Cancellation Axiom. An
{\em Abelian structure} on $\cal A$ is a subclass
$\cal P \subseteq \mor{\cal A}$
whose elements are called to be {\em proper maps}. A short exact sequence
whose maps are proper is a {\em proper s.e.s}.

The class $\cal P$ is to be subject to the following
axioms:
\begin{itemize}
\item[(P0)] $1_A\in \cal P$ for all $A\in \ob \cal A$.
\item[(P1)] If $(f:B\to C)\in \cal P$, $(g:A\to B)\in \cal P$ and
$g$ is an epimorphism then $f\circ g\in{\cal P}$; dually, if $f\in \cal P$,
$g\in\cal P$ and $f$ is a monomorphism then $f\circ g\in \cal P$.
\item[(P2)] If $f\circ g\in \cal P$ is a monomorphism then $g$ is proper;
dually, if $f\circ g\in \cal P$ is an epimorphism then $f$ is proper.
\item[(P3)] For every proper map $f:B\to D$ there are proper s.e.s
$0 \to A \to B \to C \to 0$ and $0 \to C \to D \to E \to 0$ such that
the following diagram
$$
\begin{CD}
0@>>> A @>>> B @>>> C @>>>0\\
@.  @.   @VVfV  @VV1_CV\\
0@<<< E @<<< D @<<< C @<<< 0
\end{CD}
$$
is commutative.
\item[(P4)] If in the commutative diagram
$$
\begin{CD}
  @.   0 @.  0@.  0\\
@.  @VVV @VVV @VVV\\
0@>>> A' @>>> A @>>> A'' @>>> 0\\
@.  @VVV @VVV @VVV\\
0@>>> B' @>>> B @>>> B'' @>>> 0\\
@.  @VVV @VVV @VVV\\
0@>>> C' @>>> C @>>> C'' @>>> 0\\
@.  @VVV @VVV @VVV\\
  @.   0 @.  0@.  0
\end{CD}
$$
all collumns and the second two rows are proper s.e.s, then the first
row is a proper s.e.s.
\end{itemize}
We remark that (P3)
implies the existence of kernels and cokernels
for all proper maps by \cite[Prop. 2.2]{hel1958}.

Let $\cal A$ be an additive category
with Cancellation Axiom, ${\Bbb C}$ a small category. We denote by
${\cal A}^{\Bbb C}$ the category of functors ${\Bbb C}\to {\cal A}$.
If $\rho:F\to G$ is a retraction in ${\cal A}^{\Bbb C}$ then
$\rho_c:F(c)\to G(c)$ are retractions for all $c\in \ob {\Bbb C}$ and have
kernels by Lemma 2.1. Consequently, $\rho:F\to G$ has a kernel. Therefore
${\cal A}^{\Bbb C}$ is the additive category with Cancellation Axiom.

The Abelian structure generalizes the proper class in an Abelian category in the
sense of \cite{mac1975}. We refer to \cite{kuz1972} and \cite{hel1958}
for other examples of Abelian structure.

By \cite[Prop. 3.5]{hel1958} an additive category
carries Abelian structures if
and only if it satisfies the Cancellation Axiom.

\begin{proposition}
Let $\cal A$ be an additive category with an Abelian structure  $\cal P$,
$\Bbb C$ a small category, ${\Bbb C}{\cal P}$ the class of all natural
transformations $\eta:F\to G$ for which $\eta_c\in{\cal P}$ for all
$c\in \ob {\Bbb C}$. Then the class ${\Bbb C}{\cal P}$ is an Abelian structure on
${\cal A}^{\Bbb C}$.
\end{proposition}
{\sc Proof.}
It is clear that ${\cal A}^{\Bbb C}$ is the category with Cancellation Axiom
by Lemma \ref{l21}. We check the axioms (P0)-(P4) for ${\Bbb C}{\cal P}$:
\begin{itemize}
\item[(P0)] Each identity $1_F:F\to F$ consists of $1_{F(c)}$,
$c\in \ob {\Bbb C}$. Hence $1_F\in{\cal P}$.
\item[(P1)] A morphism $\eta:F \to G$  in ${\cal A}^{\Bbb C}$ is a
monomorphism if and only if $\ker \eta=0$. Hence $\eta $ is a monomorphism
if and only if $\eta_c$ are monomorphisms for all $c\in\ob{\Bbb C}$;
dually for epimorphisms.
\item[(P2)] Analogously.
\item[(P3)] For every $f\in {\cal P}$ we choose a kernel
($\ker f$,$k_f$) of $f$.
By definition of kernels \cite[Ch.IX]{mac1975} for every commutative diagram
\begin{equation}
\begin{CD}
A @>f>> B\\
@V\alpha VV  @VV\beta V \\
A' @>f'>> B'
\end{CD}
\end{equation}
with $f,f'\in {\cal P}$ and $\alpha, \beta\in \mor{\cal A}$ there exists
unique $\ker f\to \ker f'$ for which the following diagram is commutative
$$
\begin{CD}
\ker f@>k_f>> A @>f>> B\\
@VVV @VV\alpha V @VV\beta V\\
\ker f' @>k_{f'}>> A' @>f'>> B'
\end{CD}
$$
Dually, we choose cokernels $(\coker f, c_f)$. Consider
the category which objects are all members of ${\cal P}$ and morphisms
are commutative diagrams (1).

We see that $\ker$ and $\coker$ are functors from this category to $\cal A$.
Hence for every natural transformation $\eta:F\to G$ with $\eta_c\in {\cal P}$,
for all $c\in \ob{\Bbb C}$, we have a diagram of natural transformations
$$
\begin{CD}
0 @>>> \ker \eta @>k_{\eta}>> F @>\varepsilon>> \coker k_{\eta} @>>> 0\\
@.  @.  @VV\eta V @VV\xi V \\
0 @<<< \coker\eta @<c_{\eta}<< G @<m<<  \ker c_{\eta} @<<< 0
\end{CD}
$$
where $\xi$ is an isomorphism of functors. By \cite[Prop.3.3.]{hel1958}
all equivalences in
$\cal A$ belong to $\cal P$, therefore $\xi\circ \varepsilon$ is
a proper epimorphism, and we obtain (P3) for ${\Bbb C}{\cal P}$.
\item[(P4)] $A$ sequence of functors is exact if and only if it is exact
at every $c\in \ob {\Bbb C}$. Thus (P4) holds.
\end{itemize}
Q.E.D.

\medskip

An object $P$ is called
$\cal P$-projective if for every proper epimorphism $\varepsilon:A\to B$
and for a morphism $\alpha :P\to B$ there is $\beta :P\to A$ such that
$\varepsilon \circ \beta= \alpha$.

Let $\cal A$ be an additive category with Cancellation Axiom,
${\cal P}\subseteq \mor {\cal A}$ an Abelian structure. Following Heller we
will say that $\cal A$ has enough $\cal P$-projectives if for each
$A\in \ob {\cal A}$ there is a $\cal P$-projective object $P(A)\in \cal A$
and a proper epimorphism $P(A) \stackrel{\pi A}\to A$.

\begin{proposition}
Let $\cal A$ be an additive category with coproducts, $\cal P$ an Abelian
structure. If $\cal A$ has enough $\cal P$-projectives and coproducts of proper
epimorphisms are proper, then ${\cal A}^{\Bbb C}$ has enough
${\Bbb C}{\cal P}$-projectives.
\end{proposition}
{\sc Proof.}
We consider the set $\ob {\Bbb C}$ as the maximal discrete subcategory of
${\Bbb C}$. Let $O:{\cal A}^{\Bbb C}\to {\cal A}^{\ob {\Bbb C}}$  be
the restriction functor, $O(F)=F|_{\ob {\Bbb C}}$. The  category $\cal A$ has
coproducts, hence there is a left adjoint functor $\Lambda$ to the functor $O$.
Up to a natural isomorphism, we have $\Lambda(D)(c)=\sum_{c_0\to c}D(c_o)$
for each $D\in {\cal A}^{\ob{\Bbb C}}$. The counit of adjunction
$\varepsilon_F:\Lambda O F\to F$ splits on each $c\in \ob {\Bbb C}$, hence
$\varepsilon_F$ is a proper epimorphism in ${\cal A}^{\Bbb C}$ with respect
to Abelian structure ${\Bbb C}{\cal P}$. Choose a  family
$P=\{P(c)\}_{c\in\ob{\Bbb C}}\in {\cal A}^{\ob{\Bbb C}}$ of
$\cal P$-projectives and a family of proper epimorphisms
$\psi=\{\psi_c:P(c)\to F(c)\}_{c\in\ob{\Bbb C}}$ and apply the functor
$\Lambda$. Then
$\Lambda(\psi)_c:\Lambda(P)(c)\to \Lambda (OF)(c)$ is a proper epimorphism
for each $c\in \ob {\Bbb C}$ as the coproduct of proper epimorphisms.

Consequently, $\varepsilon_F\circ(\Lambda(\psi)):\Lambda(P)\to F$
is a proper epimorphism in ${\cal A}^{\Bbb C}$.

The functor ${\cal A}^{\Bbb C}(\Lambda (P),-)$
is isomorphic to the functor ${\cal A}^{\ob {\Bbb C}}(P,O(-))$
which is exact on all proper s.e.s in
${\cal A}^{\Bbb C}$. Thus, $\Lambda (P)$  is ${\Bbb C}{\cal P}$-projective
and there is an proper epimorphism $\Lambda (P)\to F$.
Q.E.D.
\newpage

\section{COHOMOLOGY OF CATEGORIES}
Let ${\Bbb C}$ be a small category, $F:{\Bbb C}\to \ab$
a functor.

Consider the sequence of groups
$$
C^n({\Bbb C}, F)=\prod_{c_0\to \dots \to c_n}F(c_n), \quad n\geq 0.
$$

Let $N_n{\Bbb C}$ be the set of all sequences of morphisms $c_0\to c_1\to \dots \to c_n$
in $\Bbb C$. Regarding each $\varphi \in \prod_{c_0\to \dots \to c_n}F(c_n)$
as a function $\varphi :N_n {\Bbb C}\to \bigcup_{c\in \ob {\Bbb C}}F(c)$
with $\varphi (c_0\to  \dots \to c_n)\in F(c_n)$, we define the
homomorphisms $d^n:C^n({\Bbb C}, F)\to C^{n+1}({\Bbb C}, F)$ by the
formulas
$$
\begin{array}{c}
(d^n\varphi)(c_0
\stackrel{\alpha_1}\to c_1 \to \dots \stackrel{\alpha_{n+1}}\to 
c_{n+1})=\\
\sum\limits_{i=0}^n(-1)^{i}\varphi
(c_0 \stackrel{\alpha_1}\to \dots \stackrel{\alpha_i}\to \hat{c}_i
\stackrel{\alpha_{i+1}}\to 
\dots \stackrel{\alpha_{n+1}}\to 
 c_{n+1})+\\
(-1)^{n+1} F(\alpha_{n+1})(\varphi(c_0 \stackrel{\alpha_1}\to 
\dots \stackrel{\alpha_n}\to c_n)).
\end{array}
$$
Here
$(c_0 \stackrel{\alpha_1}\to \dots \stackrel{\alpha_i}\to \hat{c}_i
\stackrel{\alpha_{i+1}}\to 
\dots \stackrel{\alpha_{n+1}}\to c_{n+1})$
is equal to
$(c_0 \stackrel{\alpha_1}\to \dots \to c_{i-1} 
\stackrel{\alpha_{i+1}\circ\alpha_i}\longrightarrow c_{i+1} \to
\dots \stackrel{\alpha_{n+1}}\to 
 c_{n+1})
$
if $0< i< n+1$, and to
$(c_1 \stackrel{\alpha_2}\to \dots \stackrel{\alpha_{n+1}}\to c_{n+1})$ if
$i=0$. It is well known that $d^{n+1}\circ d^n=0$
for all $n\geq 0$. We have obtained the complex of Abelian groups and
homomorphisms
$$
0 \to C^0({\Bbb C}, F)\stackrel{d^0}\to \dots \to C^n({\Bbb C},F)
\stackrel{d^n}\to C^{n+1}({\Bbb C},F) \to \dots
$$
which is denoted by $C^*({\Bbb C}, F)$. The cohomologies
$H^n(C^*({\Bbb C}, F))=\ker d^n/\im d^{n-1}$
are isomorphic to Abelian groups $\lim\nolimits_{\Bbb C}^n F$
(see \cite{oli1992}, for example)
where $\lim\nolimits_{\Bbb C}^n:\ab^{\Bbb C}\to \ab$ are the $n$-th right derived
functors of the limit $\lim\nolimits_{\Bbb C}:\ab^{\Bbb C}\to \ab$.

In particular, if ${\Bbb C}$ has the initial object then $H^n(C^*({\Bbb C},F))=0$
for $n>0$  and $H^0(C^*({\Bbb C}, F))\cong \lim\nolimits_{\Bbb C}F$.

For each $c\in \ob  {\Bbb C}$ we denote by $c/{\Bbb C}$ the
{\em comma-category} in the sense of \cite{mac1971},
its
objects are pairs $(a,\alpha)$ of $a\in \ob  {\Bbb C}$ and $\alpha \in
{\Bbb C}(c,a)$, and for each pair of objects
$(a, \alpha \in {\Bbb C}(c,a))$ and $(b,\beta\in {\Bbb C}(c,b))$
the set of morphisms from $(a, \alpha)$ to $(b, \beta)$ consists of
$\gamma\in {\Bbb C}(a,b)$ for which  the following diagrams
\begin{equation}
\begin{CD}
c @>\alpha >> a\\
@VV=V  @VV\gamma V \\
c @>>\beta> b
\end{CD}
\end{equation}
are commutative.

Let $Q_c:c/{\Bbb C}\to {\Bbb C}$ be the functor which carries the above
diagram (2)  to $\gamma :a\to b$. The category $c/{\Bbb C}$ contains
the initial object $(c,1_c)$. Therefore,
the functor $\lim\nolimits_{c/{\Bbb C}}$ is exact and 
$H^n(C^*({\Bbb C},GQ_c))=0$ for every functor  $F:{\Bbb C}\to \ab$ and
$n>0$.

Let ${\Bbb C}$ be a small category. For each $a\in \ob {\Bbb C}$ we identify
the morphism $1_a$ with the object a, so $\ob {\Bbb C}\subseteq \mor {\Bbb C}$.
Objects of the {\em factorization category} \cite{bau1985} ${\Bbb C}'$
are all morphisms of ${\Bbb C}$, and the set of morphisms ${\Bbb C}'(f,g)$
for any $f,g\in \mor {\Bbb C}$ consists of pairs $(\alpha, \beta)$
of morphisms in $\Bbb C$ for which the diagram
$$
\begin{CD}
b@>\beta>> d\\
@AfAA  @AAgA\\
a@<<\alpha< c
\end{CD}
$$
is commutative; we denote the morphisms by $(\alpha, \beta):f\to g$.
The composition is defined  as
$(\alpha_2, \beta_2)\circ (\alpha_1,\beta_1)=(\alpha_1\circ\alpha_2,\beta_2\circ\beta_1)$.
Functors $F:€{\Bbb C}'\to \ab$ are called {\em natural systems}.

Let $\Bbb C$ be a small category, $F:{\Bbb C}'\to \ab$ a natural system. For
every $n\geq 0$ we consider an Abelian group
$$
K^n({\Bbb C},F)=
\prod_{c_0\stackrel{\alpha_1}\to c_1\stackrel{\alpha_2}\to c_2 \to \dots
 \stackrel{\alpha_n}\to c_n}
F(\alpha_n\circ\alpha_{n-1}\circ \dots \circ\alpha_1)
$$
We regard elements of  $K^n({\Bbb C},F)$ as maps
$\varphi :N_n{\Bbb C}\to \bigcup_{g\in\mor {\Bbb C}}F(g)$ with
$\varphi(\alpha_1, \dots , \alpha_n)\in F(\alpha_n\circ\alpha_{n-1}\circ \cdots\circ \alpha_1)$
and define homomorphisms $d^n:K^n({\Bbb C},F)\to K^{n+1}({\Bbb C}, F)$
by the formula
$$
\begin{array}{c}
(d^n\varphi)(\alpha_1, \dots, \alpha_{n+1})=F(\alpha_1,1)\varphi
(\alpha_2,\dots,\alpha_{n+1})+\\
\sum\limits_{i=1}^n(-1)^{i}\varphi(\alpha_1,\dots
\alpha_{i-1},\alpha_{i+1}\circ\alpha_{i},\alpha_{i+2},\dots,\alpha_{n+1})
+\\
(-1)^{n+1}F(1,\alpha_{n+1})\varphi(\alpha_1,\dots,\alpha_n).
\end{array}
$$
Cohomology groups $H^n(K^*({\Bbb C}, F))=\ker d^n/\im d^{n-1}$ are called
$n$-th cohomology groups $H^n({\Bbb C}, F)$ of ${\Bbb C}$ with coefficients
in the natural system $F$.

For any $\alpha \in \mor{\Bbb C}$ we denote by $\dom\alpha$
and $\cod\alpha$ the domain and the codomain.
We have the functor
$(\dom, \cod): {\Bbb C}'\rightarrow {\Bbb C}^{op}\times{\Bbb C}$
which acts as
$f \mapsto (\dom f, \cod f)$ on objects, and
$((\alpha, \beta): f \rightarrow g)$ $\mapsto (\alpha, \beta)$
on morphisms of ${\Bbb C}'$.
Baues and Wirsching \cite{bau1985} proved that
$H^n({\Bbb C},F)\cong\lim\nolimits_{\Bbb C'}^nF$ for all $n\geq 0$.
Functors $B:{\Bbb C}^{op}\times {\Bbb C}\to \ab$ are called
{\em bimodules} over $\Bbb C$.
For each bimodule $B$ over $\Bbb C$ we
denote by $\{B(\dom\alpha, \cod\alpha)\}$ the natural system
$B\circ (\dom,\cod):{\Bbb C}'\to \ab$. The {\em Hochschild-Mitchell
cohomologies $H^n({\Bbb C},B)$ of the category} $\Bbb C$ with coefficients
in the bimodule $B$ are the Baues-Wirsching cohomologies of
$\Bbb C$ with coefficients in the natural system $B\circ(\dom,\cod)$.
It is easy to see that $H^n({\Bbb C},B)$  isomorphic to
$\lim\nolimits_{\Bbb C'}^n\{B(\dom\alpha,\cod\alpha)\}$.

\section{EXTENSION GROUPS}
Let $\cal A$ be an additive category with an Abelian structure $\cal P$,
and with enough $\cal P$-projectives. Then for every $A\in \ob {\cal A}$
we have some object $P(A)$ and some proper epimorphism
$\pi A:P(A)\to A$. We denote by $\omega A:\Omega(A)\to P(A)$ a kernel of $\pi A$.

Let $P_*(A)$ be the exact sequence of $\cal P$-projective objects and the
proper morphisms which is obtained by sticking of sequences
$$                     
\begin{array}{ccccccccc}
0 & \to & \Omega(A)   & \stackrel{\omega A}\to & P(A)
      & \stackrel{\pi A}\to & A & \to & 0\\
0 & \to & \Omega^2(A) & \stackrel{\omega\Omega(A)}\to & P(\Omega(A))
      & \stackrel{\pi\Omega(A)}\to& \Omega(A) & \to & 0\\
&&&& \cdots\\
0 & \to & \Omega^{k+1}(A) & \stackrel{\omega \Omega^k(A)}\to & P(\Omega^k(A))
      & \stackrel{\pi \Omega^k(A)}\to & \Omega^k(A) & \to & 0\\
&&&& \cdots\\
\end{array}
$$
That is $P_*(A)$ consists of morphisms and objects:
$$
\begin{array}{c}
0 \leftarrow P(A) \stackrel{\omega A\circ \pi \Omega (A)}\longleftarrow
 P(\Omega(A))
\stackrel{\omega\Omega(A)\circ \pi\Omega^2(A)}\longleftarrow P(\Omega^2(A))
\leftarrow
\dots\\
\leftarrow P(\Omega^k(A))
\stackrel{\omega \Omega^k(A)\circ \pi \Omega^{k+1}(A)}\longleftarrow
P(\Omega^{k+1}(A)) \leftarrow \dots
\end{array}
$$
Heller define  $\ext_{\cal P}^n(A,B)$ as the group
$\hm_{-n}(P_*(A), P_*(B))$ of homotopical classes of homogenous maps.
It follows from \cite[Prop. 11.6]{hel1958} that
$\ext_{\cal P}^n(A,B)\cong H^n({\cal A}(P_*(A),B))$, $\forall n\geq 0$.

Denote by $L:\set\to \ab$ the functor from the category of sets to the category
of abelian groups which assigns to each set $E$ the free Abelian group generated by $E$
and to each map $f:E_1\to E_2$ the homomorphism extending $f$.

\begin{lemma}\label{l41}
Let ${\Bbb C}$ be a small category, $\cal A$ an additive category with coproducts,
$F:{\Bbb C}\to {\cal A}$ a functor, $D=\{D_c\}_{c\in \ob {\Bbb C}}$
a family of objects $D_c\in \ob {\cal A}$. Then
$\lim\nolimits_{{\Bbb C}'}^n\{{\cal A}(\Lambda D(\dom\alpha), F(\cod\alpha))\}=0$,
$\forall n>0$,
where $\Lambda:{\cal A}^{\ob{\Bbb C}}\to {\cal A}^{\Bbb C}$ is the left
adjoint
to the restriction functor $O:{\cal A}^{\Bbb C}\to {\cal A}^{\ob{\Bbb C}}$.
\end{lemma}
{\sc Proof.}
For every $A\in \ob{\cal A}$ and a set $E$ we denote by $\Sigma_E A$ the
coproduct of the family $\{A_e\}_{e\in E}$ of objects $A_e=A$. Let
$\nu_e:A\to \Sigma_E A$ be the canonical injections.
Consider the isomorphism
$$
\omega(E,A,B):A(\Sigma_E A,B)\to \ab (LE,{\cal A}(A,B))
$$
which assigns to every morphism $\varphi:\Sigma_E A\to B$ the homomorphism
$\omega(E,A,B):LE\to {\cal A}(A,B)$ acting an $e\in E$ as
$e\to \varphi\circ \nu_e\in {\cal A}(A,B)$.

Then we keep an arbitrary $c\in \ob {\Bbb C}$.
Consider the sets $E=h^c(a)={\Bbb C}(c,a)$ with $a\in \ob {\Bbb C}$.
Let $Lh^c$ be the composition of $L$ and $h^c$. Then there is an isomorphism
$$
{\cal A}(\Sigma_{{\Bbb C}(c,a)}A, B)\to \ab (Lh^c(a), h^{A}(B)),
$$
which is natural in $a\in {\Bbb C}$ and $A,B\in {\cal A}$. Let
$\Lambda^c:{\cal A}\to {\cal A}^{\Bbb C}$ be the left adjoint to
the functor ${\cal A}^{\Bbb C}\to {\cal A}$ acting as $F\to F(c)$
for $F\in \ob({\cal A}^{\Bbb C})$ and $\eta\to \eta_c$ for
$\eta\in \mor ({\cal A}^{\Bbb C})$. It is known that
$(\Lambda^cA)(a)=\Sigma_{{\Bbb C}(c,a)} A$ for all $a\in\ob {\Bbb C}$.
For each functor $F:{\Bbb C}\to {\cal A}$ we have an isomorphism of bifunctors
$$
{\cal A}(\Lambda^cA(-),F(=))\cong\ab(Lh^c(-),(h^{A}\circ F)(=)).
$$
It leads to an isomorphism
\begin{equation}\label{e3}
{\cal A}(\Lambda^{c}A(\dom\alpha),F(\cod\alpha))\cong
\ab(Lh^c(\dom\alpha),(h^{A}\circ F)(\cod\alpha))\quad 
\end{equation}
of the natural systems on ${\Bbb C}$ with fixed $c$.
For every family $D=\{D_c\}_{c\in\ob{\Bbb C}}$ of objects in $\cal A$ there
exists an isomorphism
$\Lambda D\cong \Sigma_{c\in{\Bbb C}}\Lambda^cD_c$. Hence
$$
\begin{array}{c}
\lim\nolimits_{{\Bbb C}'}^n\{{\cal A}(\Lambda D(\dom\alpha),F(\cod\alpha))\} \cong\\
\\
\lim\nolimits_{{\Bbb C}'}^n\{\prod_{c\in {\Bbb C}}{\cal A}
(\Lambda^cD_c(\dom\alpha),F(\cod\alpha))\}\cong\\
\\
\prod _{c\in {\Bbb C}}\lim\nolimits_{{\Bbb C}'}^n\{{\cal A}(\Lambda^cD_c
(\dom\alpha),F(\cod\alpha))\}
\end{array}
$$
\\
Thus, by the isomorphism (\ref{e3}) it suffices to prove that
$$
\lim\nolimits_{{\Bbb C}'}^n\{\ab(Lh^c(\dom\alpha),(h^{A}\circ F)
(\cod\alpha))\}=0,
$$
for every $A\in  \ob {\cal A}$ and $n > 0$. We will prove that
$$
\lim\nolimits_{{\Bbb C}'}^n\{\ab (Lh^c(\dom\alpha),G(\cod\alpha))\}=0,
$$
for every $c\in\ob {\Bbb C}$ and $G\in \ab^{\Bbb C}$.

For that purpose we consider the natural system
$$
M=\{\ab(Lh^c(\dom\alpha),G(\cod\alpha))\}
$$
and the complex $K^*({\Bbb C}, M)$. Then
$$
K^n({\Bbb C}, M)
=\prod_{c_0\to \dots\to c_n}\ab(Lh^c(c_0),G(c_n))
\cong\prod_{c_0\to \dots\to c_n}\prod_{c\to c_0}G(c_n).
$$
We recall that $C^n(c/{\Bbb C},GQ_c)$ consists of functions
$$
g:N_n(c/{\Bbb C})\to   \bigcup_{\alpha\in \ob(c/{\Bbb C})}GQ_c(\alpha)
$$
with $g(c\to c_0\to \dots \to c_n)\in GQ_c(c\to c_n)$; the
homomorphisms $d^n:C^n(c/{\Bbb C},GQ_c)\to C^{n+1}(c/{\Bbb C},GQ_c)$
act as
$$
\begin{array}{c}
(d^ng)(c\to c_0 \to \dots \to c_{n+1})=
\sum\limits_{i=0}^n(-1)^{i} g(c\to c_0\to \dots 
\to \hat{c}_{i}\to \dots \\
\to c_{n+1})+
(-1)^{n+1}G(c_n\to c_{n+1})g(c\to c_0\to \dots \to c_n).
\end{array}
$$

We will prove that the complex $K^*({\Bbb C},M)$
for
$$
M = \{\ab (Lh^c(\dom\alpha),G(\cod\alpha))\}
$$
is isomorphic to $C^*(c/{\Bbb C},GQ_c)$.

It is clear that $K^n({\Bbb C},M) \cong C^n(c/{\Bbb C},GQ_c)$.
The homomorphism $d^n :K^n({\Bbb C}, M)\to K^{n+1}({\Bbb C}, M)$
has the following action on functions $f$ for which $f(c_0\to \dots \to
c_n)(c\to c_0)\in G(c_n)$:
$$
\begin{array}{rl}
d^nf & (c_0\to c_1 \to \dots \to c_{n+1})(c\to c_0)=\\
     & f(c_1 \to \dots \to c_{n+1})\circ Lh^c(c_0 \to c_1)(c\to c_0)+\\
     & \sum_{i=1}^n(-1)^{i}f(c_0\to \dots\to \hat{c_{i}}\to \dots \to c_{n+1})(c\to c_0)+\\
     & (-1)^{n+1}F(c_n\to c_{n+1})f(c_0\to \dots \to c_n)(c\to c_0).
\end{array}
$$
We let $\tilde{f}(c\to c_0\to \dots \to c_n)=f(c_0\to \dots\to c_n)(c\to c_0)$.
We check that the correspondence $f\to \tilde{f}$
is a morphism of complexes:
$$
\begin{array}{c}
d^n\tilde{f}(c\to c_0 \to \dots \to c_{n+1})=\tilde{f}(c\to \hat{c_{0}}\to c_1\to \dots \to c_{n+1}) \\
+(-1)^{i}\tilde{f}(c\to c_0\to \dots \to \hat{c_{i}}\to \dots \to c_{n+1})+\\
(-1)^{n+1}F(c_{n}\to c_{n+1})\tilde{f}(c\to c_0\to \dots \to c_n).
\end{array}
$$
Thus, the map $f\to \tilde{f}$ is an isomorphism of complexes
$K^*({\Bbb C},M)$ and $C^*(c/{\Bbb C},GQ_c)$.
Consequently, the $n$-th cohomology groups of $K^*({\Bbb C}, M)$ are zeros
for $n>0$. Hence $\lim\nolimits_{{\Bbb C}'}^n\{\ab(Lh^c(\dom\alpha),G(\cod\alpha))\}=0$
for all $n>0$. The isomorphism (3)  and preservation of products by
$\lim\nolimits_{{\Bbb C}'}^n$ finish the proof.
\begin{theorem}
Let ${\cal A}$ be an additive category with coproducts, ${\cal P}$ a nonempty
Abelian structure. If the coproduct of every family of proper epimorphisms
is proper, and $\cal A$ has enough $\cal P$-projectives, then for each
small category $\Bbb C$  and functors $F,G:{\Bbb C}\to{\cal A}$
there exists a first quadrant spectral sequence with
$$
E_2^{p,q}=\lim\nolimits_{{\Bbb C}'}^p
\{\ext_{\cal P}^q(F(\dom\alpha),G(\cod\alpha))\}
\Longrightarrow\ext_{{\Bbb C}{\cal P}}^{p+q}(F,G).
$$
\end{theorem}
{\sc Proof.}
We will build an exact sequence
$$
0 \leftarrow F \leftarrow F_0 \leftarrow F_1 \leftarrow \dots
$$
of proper morphisms with ${\Bbb C}{\cal P}$-projective $F_n$ for all $n\geq 0$.
Recall that $O:{\cal A}^{\Bbb C}\to {\cal A}^{\ob{\Bbb C}}$
is the restriction functor,
and $\Lambda$ the left adjoint to $O$.
The counit of adjunction $\varepsilon_F:\Lambda O F\to F$ is a retraction on
each $c\in \ob {\Bbb C}$. We choose a family $P=\{P(c)\}_{c\in \ob {\Bbb C}}$
of $\cal P$-projectives and a family of proper epimorphisms
$\psi=\{\psi_c:P(c)\to F(c)\}_{c\in\ob{\Bbb C}}$
and apply the functor $\Lambda$. Then $\Lambda(\psi):\Lambda(P)\to \Lambda O F$
is the proper epimorphism as a coproduct of proper epimorphisms.
We let $F_0=\Lambda(P)$. Let $K$ be a kernel of the morphism
$\varepsilon\circ \Lambda (\psi):F_0\to F$. We apply above building to $K$
instead of $F$ and obtain some functor $K_0$. We let $F_1=K_0$.
Then the inclusion $K\to F_0$  is a proper monomorphism in ${\cal A}^{\Bbb C}$.
Consequently, the composition $F_1\to K \to F_0$ is a proper morphism.
We denote it by $d_0$. By induction we build the members
$F_2,F_3,\dots $ and morphisms $d_n:F_{n+1}\to F_n$.
Now we consider the complex
$$
\{K^n(\alpha)\}=\{{\cal A}(F_n(\dom\alpha),G(\cod\alpha))\}
$$
in the category $\ab^{{\Bbb C}'}$ of natural systems on $\Bbb C$. By
Grothendieck \cite{gro1957} there are two spectral sequences concerned with hyper cohomologies
of the functor $\lim\nolimits_{{\Bbb C}'}$ with respect to the complex $K^*$
$$
H^p(\lim\nolimits_{{\Bbb C}'}^q\{K^*(\alpha)\})\Longrightarrow H^n,\quad
\lim\nolimits_{{\Bbb C}'}^p \{H^qK^*(\alpha)\}\Longrightarrow H^n.
$$
For each $n\geq 0$ there exists a family $A = \{A_c\}_{c\in Ob{\Bbb C}}$
of objects of $\cal A$
 such that $F_n=\Lambda A$. Applying
Lemma \ref{l41} we obtain $\lim\nolimits_{{\Bbb C}'}^q\{K^n(\alpha)\}=0$
for all $q>0$, and $n\geq 0$. Thus, the first spectral sequence degenerates
and second gives the looking spectral sequence.

The {\em preabelian category} is the additive category 
with kernels and cokernels.
A morphism $f: A \rightarrow B$ in the preabelian category is called
{\em strict} if the associated morphism $Coim\, f \rightarrow Im\, f$
is the isomorphism. Kuz'minov and Cherevikin proved that any 
quasiabelian category has the Abelian structure in the sense of Heller
where the class of proper morphisms consists of all strict morphisms.
In particular the category of locally convex spaces and coninuous linear
maps is quasiabelian. Palamadov proved that this category has enough
injectives. The proper monomorphisms in this category are all the kernels.
This category has the infinite products. The product of kernels is the
kernel. Hence, the opposite to the conditions of Theorem 4.2 are satisfied.
\begin{corollary}
Let ${\Bbb C}$ be a small category and $\cal A$ the category 
of locally convex spaces and continuous linear maps.
Then for every diagrams $F,~G:~ {\Bbb C}\rightarrow {\cal A}$ 
 there exists the spectral 
sequence with 
$$
E_2^{p,q}=\lim\nolimits_{{\Bbb C}'}^p
\{\ext_{\cal P}^q(F(\dom\alpha),G(\cod\alpha))\}
\Longrightarrow\ext_{{\Bbb C}{\cal P}}^{p+q}(F,G).
$$
where $\cal P$ is the class of all strict morphisms.
\end{corollary}

\end{document}